\documentclass[conference, final, twocolumn,10pt]{IEEEtran}
\usepackage{cite}
\usepackage{graphicx}
\usepackage{mathrsfs}
\usepackage{amsmath}
\usepackage{amsfonts}
\usepackage{amssymb}

\newtheorem{theorem}{{Theorem}}
\newtheorem{lemma}{{Lemma}}
\newtheorem{corollary}{{Corollary}}
\newtheorem{definition}{{Definition}}
\newtheorem{proposition}{{Proposition}}

\newcommand{\mb}{\mathbf}
\newcommand{\qed}{\hspace*{\fill} $\Box$ \\}

\allowdisplaybreaks

\begin{document}

\title{Analytical Lower Bounds on the Critical Density in Continuum Percolation}

\author{\authorblockN{Zhenning Kong and Edmund M. Yeh}
\authorblockA{Department of Electrical Engineering\\
Yale University\\
New Haven, CT 06520, USA \\
Email: \{zhenning.kong, edmund.yeh\}@yale.edu}}

\maketitle

\begin{abstract}
Percolation theory has become a useful tool for the analysis of
large-scale wireless networks. We investigate the fundamental
problem of characterizing the critical density $\lambda_c^{(d)}$ for
$d$-dimensional Poisson random geometric graphs in continuum
percolation theory. By using a probabilistic analysis which
incorporates the clustering effect in random geometric graphs, we
develop a new class of analytical lower bounds for the critical
density~$\lambda_c^{(d)}$ in $d$-dimensional Poisson random
geometric graphs. The lower bounds are the tightest known to date.
In particular, for the two-dimensional case, the analytical lower
bound is improved to~$\lambda^{(2)}_c \geq 0.7698...$. For the
three-dimensional case, we obtain $\lambda^{(3)}_c \geq 0.4494...$.
\end{abstract}

\section{Introduction}
Recently, percolation theory has become a useful tool for the analysis of large-scale wireless
networks \cite{BoNrFrMe03,FrBoCoBrMe05,DoFrTh05,DoBaTh05}. A percolation process resides in a
random graph structure, where nodes or links are randomly designated as either ``occupied'' or
``unoccupied.'' When the graph structure resides in continuous space, the resulting model is
described by continuum percolation \cite{Gi61,Ha85,MeRo96,Pe03,Gr99,BoRi06}.  A major focus of
continuum percolation theory is the $d$-dimensional random geometric graph induced by a Poisson
point process with constant density $\lambda$. A fundamental result for continuum percolation
concerns a phase transition effect whereby the macroscopic behavior of the system is very
different for densities below and above some critical value $\lambda_c^{(d)}$. For $\lambda <
\lambda_c^{(d)}$ (subcritical), the component containing the origin contains a finite number of
points almost surely. For $\lambda> \lambda_c^{(d)}$ (supercritical), the component containing the
origin contains an infinite number of points with a positive
probability~\cite{MeRo96,Pe03,Gr99,BoRi06}.

Naturally, the characterization of the critical density $\lambda_c^{(d)}$ is a central problem in
continuum percolation theory. Unfortunately, the exact value of $\lambda_c^{(d)}$ is very
difficult to find. For two-dimensional random geometric graphs, simulation studies show that
$\lambda_c^{(2)}\approx 1.44$ \cite{QuToZi00}, while the best analytical bounds obtained thus far
are $0.696 < \lambda_c^{(2)} < 3.372$ \cite{Ha85, MeRo96}. Recently, in \cite{BaBoWa05}, the
authors reduce the problem of characterizing $\lambda_c^{(2)}$ to evaluating numerical integrals
using a mapping between continuum percolation and dependent bond percolation on lattices. By Monte
Carlo methods, they obtain numerical bounds $1.435 < \lambda_c^{(2)} < 1.437$ with confidence
$99.99\%$. Unfortunately, the bounds obtained in \cite{BaBoWa05} are not in closed-form and cannot
be generalized to higher dimensional cases.

In this paper, we give a new mathematical characterization of the critical density
$\lambda_c^{(d)}$ for Poisson random geometric graphs in $d$-dimensional Euclidean space, where $d
\geq 2$. We develop an analytical technique based on probabilistic methods~\cite{JaLuRu00} and the
{\em clustering effect} in random geometric graphs. This analysis yields a new class of lower
bounds:
\begin{equation}
\lambda_c^{(d)}\geq \frac{1}{V^{(d)} \left(1-C_t^{(d)}\right)}
\end{equation} where $V^{(d)}$ is the volume of a $d$-dimensional
unit sphere and $C_t^{(d)}$ is the $t$-th order {\em cluster
coefficient} ($t \geq 3$) for $d$-dimensional Poisson random
geometric graphs, which we will define later. This class of
\emph{analytical} lower bounds are the tightest known to date. In
particular, by evaluating $C_3^{(2)}$ in closed-form, the analytical
lower bound for two-dimensional Poisson random geometric graphs is
improved to~$\lambda^{(2)}_c \geq 0.7698...$. For three-dimensional
Poisson random geometric graphs, we obtain the analytical lower
bound $\lambda^{(3)}_c \geq 0.4494...$. By successively evaluating
$C_t^{(d)}$ for $t\geq 4$, we can obtain tighter lower bounds on
$\lambda_c^{(d)}$.

\section{Random Geometric Graphs}

In wireless networks, a communication link exists between two nodes
if the distance between them is sufficiently small, so that the
received power is large enough for successful decoding. A
mathematical model for this scenario is as follows. Let $\|\cdot\|$
be the Euclidean norm, and $f$ be some probability density function
(p.d.f.) on $\mathbb{R}^d$. Let ${\mb X}_1, {\mb X}_2, ..., {\mb
X}_n$ be independent and identically distributed (i.i.d.)
$d$-dimensional random variables with common density $f$, where
${\mb X}_i$ denotes the random location of node $i$ in
$\mathbb{R}^d$. The \emph{ensemble} of all the graphs with
undirected links connecting all those pairs $\{{\mb x}_i, {\mb
x}_j\}$ with $\|{\mb x}_i- {\mb x}_j \|\leq r, r>0,$ is called a
\emph{random geometric graph}~\cite{Pe03}. In the following, we
focus on random geometric graphs with ${\mb X}_1, {\mb X}_2, ...,
{\mb X}_n$ distributed i.i.d. according to a uniform distribution
over a given $d$-dimensional box ${\cal
A}=[0,\sqrt[d]{n/\lambda}]^d$ with $\lambda>0$. We denote such
graphs by $G({\cal X}^{(d)}_n; r)$.

Let $A=|{\cal A}|$ be the $d$-dimensional Lebesgue measure (or volume) of ${\cal A}$. An event is
said to be asymptotic almost sure (abbreviated a.a.s.) if it occurs with a probability converging
to 1 as $n \rightarrow \infty$. Consider a graph $G=(V,E)$, where $V$ and $E$ denote the set of
nodes and links respectively. Given $u,v\in V$, we say $u$ and $v$ are \emph{adjacent} if there
exists an link between $u$ and $v$, i.e., $(u,v)\in E$. In this case, we also say that $u$ and $v$
are \emph{neighbors}.

\subsection{Preliminaries}

In $G({\cal X}^{(d)}_n; r)$, let the location ${\mb x}_j$ of node
$j$ be given. A second node $i$ is randomly placed in ${\cal A}$
according to the uniform distribution $f$. There exists a link
between these two nodes if and only if node $i$ lies within a sphere
of radius $r$ around ${\mb x}_j$. Let this spherical region be
denoted as ${\cal A}({\mb x}_j)$, then the probability for the
existence of a link between $i$ and $j$ is given by
\begin{equation}\label{P-link}
P_{link}(\mathbf{x}_j)=\int_{{\cal A}(\mathbf{x}_j)}f(\mathbf{y})d\mathbf{y}.
\end{equation}
Since the underlying distribution is uniform, the probability
$P_{link}(\mathbf{x}_j)$ depends only on the volume of the
intersection between ${\cal A}$ and the node coverage volume ${\cal
A}(\mathbf{x}_j)$. Throughout this paper, we ignore border
effects\footnote{More rigorously, we may use a torus instead of a
box for ${\cal A}$. Asymptotically, as $n\rightarrow \infty$ and $A
\rightarrow \infty$ with $n/A = \lambda$ fixed, uniform random
geometric graphs on the torus are the same as those in the box.}. As
a consequence, $P_{link}({\mb x}_j)$ is independent of ${\mb x}_j$,
and thus independent among all the links:
\begin{equation}\label{Plink}
P_{link} = \frac{V^{(d)} r^d}{A} = \frac{\lambda V^{(d)} r^d}{n},
\end{equation}
where $V^{(d)}$ is the volume of a $d$-dimensional unit sphere
$V^{(d)}=\frac{\pi^{d/2}}{\Gamma(\frac{d+2}{2})}$ and $\Gamma(\cdot)$ is the Gamma function
$\Gamma(x)=\int_0^{\infty}t^{x-1}e^{-t}dt$.

It follows that in $G({\cal X}^{(d)}_n; r)$, the probability that the given node $j$ has degree
$k, 0 \leq k \leq n-1,$ is given by the binomial distribution:
\begin{equation}\label{p-kx-binomial}
p_{k}=\binom{n-1}{k}P_{link}^k[1-P_{link}]^{n-1-k}.
\end{equation}
Thus, the mean degree for each node is
\begin{equation}\label{E-k-uniform}
\mu=E[k]= (n-1)P_{link} = \frac{(n-1)\lambda V^{(d)} r^d}{n}.
\end{equation}

Note that as $n$ and $A$ both become large but the ratio $n/A=\lambda$ is kept constant, each node
has an approximately Poisson degree distribution \cite{Pe03,Fe57} with an expected degree
\begin{equation}\label{mu-lambda}
\mu=\lim_{n\rightarrow \infty}\frac{(n-1)\lambda V^{(d)} r^d}{n}=\lambda V^{(d)}r^d.
\end{equation}
As $n\rightarrow \infty$ and $A \rightarrow \infty$ with $n/A=\lambda$ fixed, $G({\cal
X}^{(d)}_n;r)$ converges in distribution to an (infinite) random geometric graph
$G(\mathcal{H}^{(d)}_{\lambda};r)$ induced by a homogeneous Poisson point process with density
$\lambda>0$.  For such graphs, we have the following lemma.

\vspace{0.1in}%
\begin{lemma}\label{Lemma-Infinite-Vertices}
Suppose $G(\mathcal{H}_{\lambda}^{(d)};r)$ is a random geometric graph induced by a homogeneous
Poisson point process with density $\lambda>0$. Let ${\cal A}^{\prime}$ be any subset of
$\mathbb{R}^d$ with $|\mathcal{A}'|<\infty$, then the subgraph $G'$ contained in ${\cal
A}^{\prime}$ has a finite number of nodes a.a.s.
\end{lemma}
\vspace{0.1in}%

\emph{Proof:} The proof is straightforward and omitted here.\qed

Due to the scaling property of random geometric graphs \cite{MeRo96,Pe03}, in the following, we
focus on $G(\mathcal{H}^{(d)}_{\lambda};1)$ and $G(\mathcal{X}^{(d)}_n;1)$.

\subsection{Cluster Coefficients}

An important characteristic of random geometric graphs is the {\em clustering effect}.  Here, if
node $i$ is close to node $j$, and node $j$ is close to node $k$, then $i$ is typically also close
to $k$.  In the following, we use the {\em cluster coefficients} to precisely characterize the
clustering property. This turns out to be the key to deriving new bounds for the critical density
in continuum percolation.

\vspace{0.1in}%
\begin{definition} Given distinct nodes $i,j,k$ in $G({\cal X}^{(d)}_n;1)$,
the \emph{cluster coefficient} $C^{(d)}$ is the conditional probability that nodes $i$ and $j$
are adjacent given that $i$ and $j$ are both adjacent to node~$k$.
\end{definition}
\vspace{0.1in}%

\begin{figure}[t]
\centering
\includegraphics[width=1.2in]{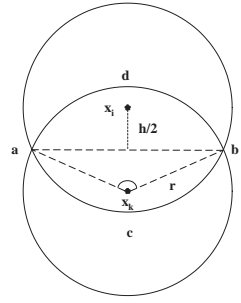}
\caption{Calculation of cluster coefficient $C^{(2)}$}
\end{figure}

The calculation of $C^{(2)}$ for two-dimensional random geometric graphs is illustrated in
Figure~1. To determine $C^{(2)}$, assume both nodes $i$ and $j$ lie within ${\cal A}({\mb x}_k)$,
then the conditional probability that nodes $i$ and $j$ are also adjacent is equal to the
probability that two randomly chosen points in a circle with radius 1 is at most distance 1 apart.
In other words, given the coordinates of $\mathbf{x}_k$ and $\mathbf{x}_i$, the probability that
there is an link between $i$ and $j$ is equal to the fraction of ${\cal A}({\mb x}_i)$ that
intersects ${\cal A}({\mb x}_k)$. By averaging ${\mb x}_i$ over all points in ${\cal A}({\mb
x}_k)$, the cluster coefficient can be found as $C^{(2)} =1-\frac{3\sqrt{3}}{4\pi}\approx
0.5865...$ \cite{DaCh02}.

The cluster coefficient $C^{(d)}$ reflects the {\em triangle effect} in random geometric graphs.
To further capture the cluster effect, we now generalize the notion of cluster coefficients for
more than three nodes.

\vspace{0.1in}%
\begin{definition} For $t\geq 3$, suppose $v_1, \ldots, v_{t-1} \in G({\cal X}^{(d)}_n; 1)$
form a single chain, i.e., they satisfy the following properties:
\begin{itemize}
\item[i)] For each $j=1,2,...,t-2$, $(v_j,v_{j+1})\in E$. \item[ii)]  For all $1\leq j,k\leq t-1$,
$(v_j,v_k)\notin E$ for $|j-k| > 1$,
\end{itemize}
where $E$ denotes the set of links in $G({\cal X}^{(d)}_n; 1)$. Then the \emph{$t$-th order
cluster coefficient $C^{(d)}_t$} is defined to be the conditional probability that a node $v_t$ is
adjacent to at least one of the nodes $v_2,...,v_{t-1},$ given that $v_t$ is adjacent to $v_1$
(averaging over all the possible positions ${\mb X_{v_2}},\ldots,{\mb X_{v_{t-1}}}$ in
$\mathbb{R}^d$ of the points $v_2$, $\ldots, v_{t-1}$ satisfying conditions (i) and (ii)).
\end{definition}
\vspace{0.1in}%

According to the above definition, $C_3^{(d)}=C^{(d)}$. To calculate $C_t^{(d)}$ is difficult in
general. However, $C^{(d)}_3$, the cluster coefficient for $d$-dimensional Poisson random
geometric graphs can be computed as \cite{DaCh02}
\begin{equation}\label{L-3d}
C^{(d)}_3=\frac{3}{\sqrt{\pi}}\frac{\Gamma(\frac{d+2}{2})}{\Gamma(\frac{d+1}{2})}\int_0^{\pi/3}\sin^d\theta
d\theta.
\end{equation}

Among higher dimensional random geometric graphs, the three-dimensional case is of practical
interest (as for sensor networks in the deep sea). Using the duplicate formula for the Gamma
function, we can derive $C_3^{(3)}$ from (\ref{L-3d}) as \cite{DaCh02}
\begin{equation}\label{C-3-3}
C_3^{(3)}=\frac{3}{2}-\frac{1}{\sqrt{\pi}}\sum_{i=1/2}^{3/2}\frac{\Gamma(i)}{\Gamma(i+\frac{1}{2})}\Big(\frac{3}{4}\Big)^{i+\frac{1}{2}}
=0.4688...
\end{equation}

Although it is difficult to obtain closed-from expressions for $C^{(d)}_t$, $t\geq 4$, we are able
to compute them by numerical integration. For example, we obtain $C^{(2)}_4 \approx 0.6012$ and
$C^{(2)}_5 \approx 0.6179$. We also note that $0 < C^{(d)}_t < 1 $, for all $t\geq 3$, since
$C^{(d)}_t$ is a (nonzero) conditional probability.

\subsection{Critical Density for Random Geometric Graphs}

Let $\mathcal{H}^{(d)}_{\lambda,\mathbf{0}}=\mathcal{H}^{(d)}_{\lambda}\cup \{\mathbf{0}\}$, i.e.,
the union of the origin and the infinite homogeneous Poisson point spatial process with density
$\lambda$ in $\mathbb{R}^d$. Note that in a random geometric graph induced by a homogeneous
Poisson point process, the choice of the origin can be arbitrary.

\vspace{0.1in}%
\begin{definition} For $G(\mathcal{H}^{(d)}_{\lambda,\mathbf{0}};1)$, the
\emph{percolation probability} $p_{\infty}(\lambda)$ is the probability that the component
containing the origin has an infinite number of nodes of the graph.
\end{definition}
\vspace{0.1in}%

\begin{definition} For $G(\mathcal{H}^{(d)}_{\lambda,\mathbf{0}};1)$, the
\emph{critical density} (continuum percolation threshold) $\lambda^{(d)}_c$ is defined as
$\lambda^{(d)}_c=\inf \{\lambda>0: p_{\infty}(\lambda)>0\}$.
\end{definition}
\vspace{0.1in}%

It is known from continuum percolation theory that if $\lambda>\lambda^{(d)}_c$, then there exists
a unique connected component containing $\Theta(n)$ nodes in $G(\mathcal{X}^{(d)}_n;1)$
a.a.s.\footnote{We say $f(n) = O(g(n))$ if there exists $n_0 > 0$ and constant $c_0$ such that
$f(n) \leq c_0g(n)~\forall n \geq n_0$.
We say $f(n)
= \Omega(g(n))$ if $g(n) = O(f(n))$. Finally, we say $f(n) =
\Theta(g(n))$ if $f(n) = O(g(n))$ and $f(n) = \Omega(g(n))$.} This
large connected component is called the {\em giant
component}~\cite{MeRo96}.

\section{New Lower Bounds on the Critical Density}

A fundamental result of continuum percolation states that
$0<\lambda^{(d)}_c<\infty$ for all $d \geq 2$. Exact values for
$\lambda^{(d)}_c$ and $p_{\infty}(\lambda)$ are not yet known. For
$d=2$, simulation studies \cite{QuToZi00} show that $\lambda^{(2)}_c
\approx 1.44$, while the best analytical bounds obtained thus far
are $0.696 < \lambda_c^{(2)} < 3.372$ \cite{Ha85, MeRo96}. Recently,
in \cite{BaBoWa05}, the authors reduce the problem of characterizing
$\lambda_c^{(2)}$ to evaluating numerical integrals, and they obtain
numerical bounds $1.435 < \lambda_c^{(2)} < 1.437$ with confidence
$99.99\%$.  Unfortunately, these bounds are not in closed form and
are restricted to the two dimensional case. In the following, we
present an analysis which combines a technique used
in~\cite{JaLuRu00} (for random graphs) and the {\em clustering
effect} in random geometric graphs to obtain a new mathematical
characterization of the critical density $\lambda^{(d)}_c$ for $d
\geq 2$. This analysis yields a new class of improved lower bounds
for~$\lambda^{(d)}_c$.  In particular, they yield the tightest
analytical lower bounds known to date.

\vspace{0.1in}%
\begin{theorem}\label{Theorem-RGG-Lower-Bound}
Let $\mu$ be the mean degree of $G(\mathcal{X}^{(d)}_n;1)$, where $d
\geq 2$. For any given integer $t\geq 3$, if
\[
\mu<\frac{1}{1-C^{(d)}_t},
\]
where $C^{(d)}_t$ is defined by Definition 2, then the largest component of
$G(\mathcal{X}^{(d)}_n;1)$ has at most $\alpha\ln n$ nodes a.a.s., where $\alpha$ is a positive
constant.
\end{theorem}
\vspace{0.1in}%

Before giving the proof, we define the diameter of a node with respect to a subgraph.
\vspace{0.1in}%
\begin{definition} Given a graph $G=(V,E)$, for any subgraph
$G^{\prime}=(V^{\prime},E^{\prime})\subseteq G$ and a node $u\in V^{\prime}$, define the
\emph{diameter of node $u$ with respect to $G^{\prime}=(V^{\prime},E^{\prime})$} as
\begin{equation}\label{Diameter}
\mbox{diam}(u,G^{\prime})\equiv \max_{v\in V'} \{ d(u,v)\},
\end{equation}
where $d(u,v)$ is the distance between $u$ and $v$, measured by the length of the shortest path
between $u$ and $v$ in terms of the number of links.
\end{definition}
\vspace{0.1in}%

Note that the diameter of graph $G$ is the maximum of the diameters of all nodes with respect to
graph $G$, i.e., $\mbox{diam}(G)=\max_{u\in G}\{\mbox{diam}(u,G)\}$. Another useful fact for the
following proof is that for a random geometric graph $G(\mathcal{X}_n;1)$, if
$\mbox{diam}(u,G^{\prime})\leq c$, then the Euclidean distance between $u$ and any node $v$ in
$G^{\prime}=(V^{\prime},E^{\prime})$ is no more than $c$, i.e., $||{\mb X}_u- {\mb X}_v||\leq c$,
for all $v\in V^{\prime}$

\vspace{0.1in} %
\emph{Proof of Theorem~\ref{Theorem-RGG-Lower-Bound}}: Let $\mu=\frac{1-\epsilon}{1-C^{(d)}_t}$,
where $0<\epsilon<1$. For simplicity, let $p$ denote the probability that there is an link between
two nodes, i.e. $p=P_{link}$ given by (\ref{Plink}), so that $(n-1)p=\mu$.

We consider an arbitrary node (with fixed label and random position) $v \in
G(\mathcal{X}^{(d)}_n;1)$ and study the following ``active-saturated'' process. For
$i=0,1,2,\ldots$, let $A_i$ denote the set of ``active'' nodes, and $S_i$ denote the set of
``saturated'' nodes, starting with $A_0=\{v\}, S_0=\emptyset$. At $(i+1)$-th step, we select an
arbitrary node $u$ from $A_i$ and update the active and saturated sets as follows:
\[
A_{i+1}=(A_i\backslash u)\cup(N_i\cap(A_i\cup S_i)^c), \;\;S_{i+1}=S_i\cup u,
\]
where $N_i$ is the set of neighbors of $u$. In other words, at each step we move a node $u$ from
the active set to the saturated set, and at the same time move to the active set all the neighbors
of $u$ which do not currently belong to the active or saturated set. In this manner, we can go
through all the nodes in $v$'s component, represented by $\Gamma_v$, until $A_i=\emptyset$.

Let $Y_{i+1}$ be the number of nodes added to $A_i$ at step $i+1$:
\begin{equation}\label{Y-i+1}
Y_{i+1}=|N_i\cap(A_i\cup S_i)^c|.
\end{equation}
Note that $|S_i|=i$ for $i\leq |\Gamma_v|$.

We say a sample graph $G_n$ of $G(\mathcal{X}^{(d)}_n;1)$ is {\em good}, if there is no component
having size strictly larger than $\frac{3-2\epsilon}{\epsilon^2}\ln n$, or for any node $v\in G_n$
with $|\Gamma_v|>\frac{3-2\epsilon}{\epsilon^2}\ln n$ and any sequence of the ``active-saturated"
steps starting at $v$, there exists a bounded $k'$ such that
\begin{equation}\label{k-prime-RGG}
\forall j\geq k^{\prime}\mbox{ and }\forall u \in A_j, \quad \mbox{diam}(u,S_j\cup u)\geq t-2,
\end{equation}
\begin{equation}\label{c-t-RGG}
c_{k'} \equiv |A_{k^{\prime}}\cup S_{k^{\prime}}|<(\ln n)^{1/3}.
\end{equation}
Note that $k'$ depends on $n$, the sample graph $G_n$, the node
$v\in G_n$ and the sequence of the ``active-saturated" steps
starting at $v$. Let $\mathcal{T}_n$ be the collection of all good
sample graphs with $n$ nodes. We will show later that with
probability 1, there exists a uniform bound $k_0 < \infty$ such that
$k' \leq k_0$ for any $n$, any $G_n \in \mathcal{T}_n$ and any $v\in
G_n \in \mathcal{T}_n$ with
$|\Gamma_v|>\frac{3-2\epsilon}{\epsilon^2}\ln n$ and any sequence of
the ``active-saturated" steps starting at $v$. We assume that this
holds for the moment.

Now given $G(\mathcal{X}^{(d)}_n;1) \in\mathcal{T}_n$, consider an arbitrary node (with fixed
label) $v \in G(\mathcal{X}^{(d)}_n;1)$, if $|\Gamma_v|>\frac{3-2\epsilon}{\epsilon^2}\ln n$, the
``active-saturated" process can sustain at least $\frac{3-2\epsilon}{\epsilon^2}\ln n$ steps.
Let\footnote{We ignore integer constraints for convenience.}
\[
k=\frac{3-2\epsilon}{\epsilon^2}\ln n-k'\geq 0.
\]
Because $|S_i|=i$, we have
\begin{equation}\label{A-S-kkprime}
|\Gamma_v|\geq k+k^{\prime} \Longleftrightarrow |A_{k+k^{\prime}}\cup S_{k+k^{\prime}}|\geq
k+k^{\prime}.
\end{equation}
Since $G(\mathcal{X}^{(d)}_n;1) \in\mathcal{T}_n$, with probability 1, there exist constants
$k_0$, and $k'\leq k_0$ satisfying condition (\ref{k-prime-RGG}) and (\ref{c-t-RGG}). By the
definition of $c_{k'}$, (\ref{A-S-kkprime}) is equivalent to
\begin{equation}\label{Xi-kkprime}
|\Gamma_v|\geq k+k^{\prime} \Longleftrightarrow \sum_{i=1+k^{\prime}}^{k+k^{\prime}}Y_i \geq
k+k^{\prime}-c_{k'}.
\end{equation}
Therefore,
\begin{equation}\label{Pr-Gammav}
\mbox{Pr}\{|\Gamma_v|\geq k+k^{\prime}\}= \mbox{Pr}\left\{\sum_{i=1+k^{\prime}}^{k+k^{\prime}}Y_i
\geq k+k^{\prime}-c_{k'}\right\}.
\end{equation}

We now bound the RHS probability. Consider the ``active-saturated"
process after $k^{\prime}$ steps. For all $ j\geq k^{\prime}$, we
have $\mbox{diam}(u,S_j\cup u)\geq t-2, \forall u \in A_j$. At each
step $i$, we move $Y_i$ nodes to the active set. Suppose at the
$(j+1)$-th step, $\sum_{i=1}^jY_i=m$, and there are $n-1-m$ nodes
remaining.  Now suppose we move node $u_{j+1}$ from the active set
to the saturated set, then $Y_{j+1}$ is the number of nodes adjacent
to $u_{j+1}$ but not in $A_j$ or $S_j$. Since all nodes in $A_j$ are
adjacent to some node in $S_j$, $Y_{j+1}$ is also the number of
nodes adjacent to $u_{j+1}$, not in $S_j$ and not adjacent to any
node in $S_j$. Since $\mbox{diam}(u_{j+1},S_j\cup u_{j+1})\geq t-2$,
there exists a sequence of nodes $w_{j+1}^1,...,w_{j+1}^{t-2}$ in
$S_j$ that forms a single chain with node $u_{j+1}$ (i.e., satisfies
condition (i)-(ii) of Definition 2). Let
$\tilde{C}_t^{(d)}(n,u_{j+1},w_{j+1}^1,...,w_{j+1}^{t_2})$ be the
conditional probability that one of the remaining $n-1-m$ nodes,
$w$, is adjacent to at least one of the nodes
$w_{j+1}^1,...,w_{j+1}^{t_2}$ given that $w$ is adjacent to
$u_{j+1}$. Then the probability that a node is adjacent to $u_{j+1}$
and not adjacent to any $w_{j+1}^i,i=1,2,...,t-2$, is
$q_n(u_{j+1})\equiv
p_n(u_{j+1})(1-\tilde{C}_t^{(d)}(n,u_{j+1},w_{j+1}^1,...,w_{j+1}^{t_2}))$,
where $p_n(u_{j+1})$ is the average probability that there is a link
between node $u_{j+1}$ and any other node. Since
$\mbox{diam}(u_{j+1},S_j\cup u_{j+1})$ may be larger than $t-2$, and
there are other geometric constraints for each of the $n-1-m$
remaining nodes, the probability of any one of the remaining $n-m-1$
nodes is adjacent to $u_{j+1}$, not in $S_j$ and not adjacent to any
node in $S_j$ is less than or equal to $q_n(u_{j+1})$.

Now $Y_{j+1}=\sum_{i=1}^{n-1-m}B_i$, where $B_i=1$ if node $i$ is
adjacent to $u_{j+1}$, not in $S_j$ and not adjacent to any node in
$S_j$, and $B_i=0$ otherwise.  Note that these $B_i$'s are not
independent. Nevertheless, by the argument above, we have
$\mbox{Pr}\{B_{l+1}=1 | (B_1, B_2,...,B_l)=(b_1,b_2,...,b_l) \}\leq
q_n(u_{j+1})$ for any $(b_1, b_2,...,b_l) \in \{0,1\}^l$ and $l \leq
n-2-m$, For $i=1,\ldots,n-1-m$, let $B^+_i$ be i.i.d.
Bernoulli($q_n(u_{j+1})$) random variables. By Proposition~1 in
Appendix~I, $Y_{j+1}=\sum_{i=1}^{n-1-m}B_i$ is stochastically upper
bounded by $Y''_{j+1} \equiv \sum_{i=1}^{n-1-m}B^+_i \sim
\mbox{Binom}(n-1-m,q_n(u_{j+1}))$, which is further stochastically
upper bounded by $Y'_{j+1} \sim \mbox{Binom}(n-1,q_n(u_{j+1}))$.
Therefore, conditional on $\sum_{i=1}^jY_i=m$, for any $m\leq n-1$,
$Y_{j+1}$ is stochastically upper bounded by a random variable
$Y'_{j+1}$ with distribution $\mbox{Binom}(n-1,q_n(u_{j+1}))$.

Using the same argument, we see that $\sum_{i=1+k'}^{k+k'}Y_i$ is stochastically upper
bounded by $\sum_{i=1+k^{\prime}}^{k+k^{\prime}}Y'_i\equiv Z_k \sim
\mbox{Binom}(k(n-1),q_n(u_m))$, where $q_n(u_m)= \sup_{1+k'\leq i\leq k+k'}q_n(u_i)$.
Thus,
\begin{equation}\label{Pr-Xi}
\mbox{Pr}\left\{\sum_{i=1+k^{\prime}}^{k+k^{\prime}}Y_i \geq k+k^{\prime}-c_{k'}\right\}\leq
\mbox{Pr}\{Z_k \geq k+k^{\prime}-c_{k'}\}.
\end{equation}
Let $\mu_n(u_m) \equiv (n-1)p_n(u_m)$ be the mean degree of node
$u_m$. Since
$E[Z_k]=k(n-1)q_n(u_m)=k\mu_n(u_m)(1-\tilde{C}_t^{(d)}(n,u_m,w_{m}^1,...,w_{m}^{t_2}))$,
\begin{eqnarray*}
& & \mbox{Pr}\{Z_k \geq k+k^{\prime}-c_{k'}\} \\
& = & \mbox{Pr}\{Z_k \geq E[Z_k]+k+k^{\prime}-c_{k'}-E[Z_k]\}\\
& = & \mbox{Pr}\{Z_k \geq E[Z_k] + k\delta_n(u_m)+k^{\prime}-c_{k'} \},
\end{eqnarray*}
where
$\delta_n(u_m)=1-\mu_n(u_m)(1-\tilde{C}_t^{(d)}(n,u_m,w_{m}^1,...,w_{m}^{t_2}))$.
Note that conditioned on $G(\mathcal{X}^{(d)}_n;1)\in
\mathcal{T}_n$, the node distribution may not be uniform.
Nevertheless, we will show that
$\mbox{Pr}\{G(\mathcal{X}^{(d)}_n;1)\in \mathcal{T}_n\} \rightarrow
1$ as $n \rightarrow \infty$. Hence the node distribution is uniform
asymptotically, and $p_n(u_{j+1})\rightarrow p$,
$\mu_n(u_{j+1})\rightarrow \mu$ and
$\tilde{C}_t^{(d)}(n,u_{j+1},w_{j+1}^1,...,w_{j+1}^{t_2})\rightarrow
C_t^{(d)}$ as $n \rightarrow \infty$. Thus $\delta_n(u_m)
\rightarrow \epsilon$ as $n \rightarrow \infty$. By
\eqref{c-t-prime}, there exists $0<n_0<\infty$, such that for $n\geq
n_0$, $|\delta_n(u_m)-\epsilon|\leq \frac{\epsilon}{2}$ and
\begin{eqnarray*}
& & k\delta_n(u_m)+k^{\prime}-c_{k'}\\
& = &\delta_n(u_m)\frac{3-2\epsilon}{\epsilon^2}\ln n +(1-\delta_n(u_m))k'-c_{k'}\\
& > & \frac{3-2\epsilon}{2\epsilon}\ln n-(\ln n)^{\frac{1}{3}}\\
& > & 0.
\end{eqnarray*}
By the Chernoff bound \cite{JaLuRu00}, for $\delta > 0$,
\begin{equation}\label{Chernoff-Bound}
\mbox{Pr}\{Z\geq E[Z]+\delta\}\leq \exp \left\{-\frac{\delta^2}{2E[Z]+2\delta/3}\right\}.
\end{equation}
Thus, for $n$ sufficiently large,
\begin{small}
\begin{equation}\label{Yk1}
\begin{array}{ll}
&\!\!\!\mbox{Pr}\{Z_k \geq k+k^{\prime}-c_{k'}\} \\
\leq &\!\!\!\exp
\left\{-\displaystyle\frac{(k\delta_n(u_m)+k^{\prime}-c_{k'})^2}{2k(1-\delta_n(u_m))
+2(k\delta_n(u_m)+k^{\prime}-c_{k'})/3}\right\}\\ = &\!\!\!\exp
\left\{-\displaystyle\frac{3}{2}\left[\frac{k^2\delta_n(u_m)^2+2k\delta_n(u_m)(k^{\prime}-c_{k'})
+(k^{\prime}-c_{k'})^2}{k(3-2\delta_n(u_m))+k^{\prime}-c_{k'}}\right]\right\}
\end{array}
\end{equation}
\end{small}
Since $k=\frac{3-2\epsilon}{\epsilon^2}\ln n-k'$, $k' \leq k_0$ and $c_{k'}<(\ln
n)^{\frac{1}{3}}$, as $n\rightarrow \infty$, the RHS of \eqref{Yk1} has the same order as
\begin{eqnarray}
& & \exp\Big\{-\frac{3}{2}\Big[\frac{(3-2\epsilon)\delta_n(u_m)^2}{\epsilon^2(3-2\delta_n(u_m))}
\ln n \nonumber\\
& & \qquad \qquad \quad + \frac{2\delta_n(u_m)}{3-2\delta_n(u_m)}\Big(k'\Big(1-\frac{\delta_n(u_m)}{2}\Big)-c_{k'}\Big)\Big]\Big\}\nonumber\\
& \leq & \exp\left\{-\frac{3\delta_n(u_m)(2-\delta_n(u_m))}{2(3-2\delta_n(u_m))}\right\}\nonumber\\
& & \exp\left\{-\frac{3(3-2\epsilon)\delta_n(u_m)^2}{2\epsilon^2(3-2\delta_n(u_m))}\ln n
    +\frac{3\delta_n(u_m)}{3-2\delta_n(u_m)}c_{k'}\right\}\label{Yk2}
\end{eqnarray}
Now choose $\gamma>0$ such that
$\frac{(\epsilon-\gamma)^2(3-2\epsilon)}{\epsilon^2[3-2(\epsilon-\gamma)]}>\frac{5}{6}$.
Because $\delta_n(u_m) \rightarrow \epsilon$ as $n \rightarrow \infty$, there exists
$0<n_1<\infty$, such that for $n\geq n_1$, $|\delta_n(u_m)-\epsilon|\leq \gamma$. Then
using $\epsilon-\gamma\leq\delta_n(u_m)\leq \epsilon+\gamma$, we can bound \eqref{Yk2} by
\begin{eqnarray}
& & c'\exp\left\{-\frac{3(3-2\epsilon)(\epsilon-\gamma)^2}{2\epsilon^2[3-2(\epsilon-\gamma)]}\ln n
+ \frac{3(\epsilon+\gamma)}{3-2(\epsilon+\gamma)}c_{k'}\right\}\nonumber\\
& \leq & c'\exp\left\{-\frac{5}{4}\ln
n+\frac{3(\epsilon+\gamma)}{3-2(\epsilon+\gamma)}c_{k'}\right\}\nonumber\\
&= & O(n^{-\frac{5}{4}}),\label{Yk3}
\end{eqnarray}
where
\[
c'=\exp\left\{-\frac{3(\epsilon-\gamma)[2-(\epsilon+\gamma)]}{2[3-2(\epsilon-\gamma)]}\right\}.
\]

By (\ref{Pr-Gammav})-(\ref{Yk3}), for any arbitrary node $v\in G(\mathcal{X}^{(d)}_n;1) \in
\mathcal{T}_n$,
\begin{equation}\label{Pr-Gammav-infinite}
\mbox{Pr}\Big\{|\Gamma_v|\geq \frac{3-2\epsilon}{\epsilon^2}\ln n\Big\}=O(n^{-\frac{5}{4}}).
\end{equation}
Set $\alpha=\frac{3-2\epsilon}{\epsilon^2}$. The probability that random geometric graph
$G(\mathcal{X}^{(d)}_n;1)$ has at least one component whose size is no smaller than
$\alpha\ln n$ is
\begin{eqnarray*}
\!\!\!& & \!\!\!\!\!\mbox{Pr}\{\exists v\in G(\mathcal{X}^{(d)}_n;1): |\Gamma_v|\geq \alpha\ln n\}\\
\!\!\!& = &\!\!\!\!\! \mbox{Pr}\{\exists v \in G(\mathcal{X}^{(d)}_n;1):|\Gamma_v|\geq \alpha\ln n
| G(\mathcal{X}^{(d)}_n;1)\in\mathcal{T}_n\}\\
\!\!\!& &\!\!\!\!\!\cdot\mbox{Pr}\{G(\mathcal{X}^{(d)}_n;1)\in\mathcal{T}_n\}\\
\!\!\!& &\!\!\!\!\! + \mbox{Pr}\{\exists v\in G(\mathcal{X}^{(d)}_n;1): |\Gamma_v|\geq \alpha\ln n
|G(\mathcal{X}^{(d)}_n;1)\notin\mathcal{T}_n\}\\
\!\!\!& &\!\!\!\!\! \cdot\mbox{Pr}\{G(\mathcal{X}^{(d)}_n;1)\notin\mathcal{T}_n\}\\
\!\!\!& \leq &\!\!\!\!\! n\mbox{Pr}\{|\Gamma_v|\geq \alpha\ln n |
G(\mathcal{X}^{(d)}_n;1)\in\mathcal{T}_n\}
\mbox{Pr}\{G(\mathcal{X}^{(d)}_n;1)\in\mathcal{T}_n\}\\
\!\!\!& &\!\!\!\!\! + 1\cdot\mbox{Pr}\{G(\mathcal{X}^{(d)}_n;1)\notin\mathcal{T}_n\}\\
\!\!\!& = &\!\!\!\!\!
O(n^{-\frac{1}{4}})+\mbox{Pr}\{G(\mathcal{X}^{(d)}_n;1)\notin\mathcal{T}_n\}.
\end{eqnarray*}

To complete the proof, we show two facts:
\begin{itemize}\item[(i)] $\mbox{Pr}\{G(\mathcal{X}^{(d)}_n;1)\in\mathcal{T}_n\} \rightarrow 1$ as $n \rightarrow
\infty$. \item[(ii)] With probability 1, there exists $k_0 < \infty$ such that $k' \leq k_0$ for
any $n$, any $G(\mathcal{X}^{(d)}_n;1)\in\mathcal{T}_n$ and any $v \in
G(\mathcal{X}^{(d)}_n;1)\in\mathcal{T}_n$, and any realization of the ``active-saturated" process
starting at $v$.\end{itemize}

To show (i), note that in $G(\mathcal{X}^{(d)}_n;1)$, if there is no component having size
strictly larger than $\alpha\ln n$, then it is good; otherwise, we prove that for any node $v\in
G(\mathcal{X}^{(d)}_n;1)$ with $|\Gamma_v|>\frac{3-2\epsilon}{\epsilon^2}\ln n$ and any
realization of the ``active-saturated" process starting at $v$, there exists a bounded $k'$ such
that~\eqref{k-prime-RGG} holds a.a.s. Suppose for some $v\in G(\mathcal{X}^{(d)}_n;1)$ with
$|\Gamma_v|>\frac{3-2\epsilon}{\epsilon^2}\ln n$ and a realization of the ``active-saturated"
process starting at $v$ such that for any $k\leq |\Gamma_v|-1$, there exists a step $j, k\leq j
\leq |\Gamma_v|-1$, and $w\in A_j$, such that $\mbox{diam}(w,S_j\cup w)< t-2$. Since
$|\Gamma_v|>\frac{3-2\epsilon}{\epsilon^2}\ln n$, as $n\rightarrow \infty$, $|\Gamma_v|\rightarrow
\infty$, the ``active-saturated" process can go on forever. Thus, $S_j$ asymptotically contains an
infinite number of nodes. Since $\mbox{diam}(w,S_j\cup w)< t-2$, all the nodes of $S_j$ lie in a
ball centered at $w$ with radius $t-2$ (by the argument immediately following Definition~5), which
occurs with probability approaching 0 as $n \rightarrow \infty$, by Lemma~1.  Now suppose $k'$
satisfying~\eqref{k-prime-RGG} exists but is unbounded as $n\rightarrow \infty$.  Then let
$\tilde{k}'$ be the smallest $k'$ satisfying~\eqref{k-prime-RGG}.  For step $j=\tilde{k}'-1$,
there exists at least one node $w\in A_j$ such that $\mbox{diam}(w,S_j\cup w)< t-2$.   Then,
arguing as before, we can show this happens with probability approaching 0 as $n \rightarrow
\infty$.

Next we show that $c_{k'}<(\ln n)^{\frac{1}{3}}$ a.a.s.  We know $c_{k'} = |A_{k^{\prime}}\cup
S_{k^{\prime}}|=\sum_{i=1}^{k'}Y_i$. Using an argument similar to that leading to~\eqref{Pr-Xi},
we see that the RHS is stochastically upper bounded by a random variable with distribution
$\mbox{Binom}(k'(n-1),p)$. Assuming fact (ii), which is shown below, the RHS is further
stochastically upper bounded by a random variable $\tilde{c}_{k'}$ with distribution
$\mbox{Binom}(k_0 n,p)$.  Applying the Chernoff bound (\ref{Chernoff-Bound}), we have, for
sufficiently large~$n$,
\begin{eqnarray}\label{c-t-prime}
\mbox{Pr}\{c_{k'}\geq(\ln n)^{\frac{1}{3}}\}
\!\!\!\!\!& \leq & \!\!\!\!\!\mbox{Pr}\{\tilde{c}_{k'}\geq(\ln n)^{\frac{1}{3}}\}\nonumber\\
\!\!\!\!\!& = &\!\!\!\!\!\mbox{Pr}\{\tilde{c}_{k'}\geq k_0\mu+(\ln n)^{\frac{1}{3}}-k_0\mu\}\nonumber\\
\!\!\!\!\!& \leq &\!\!\!\!\! \exp\left\{-\frac{[(\ln n)^{\frac{1}{3}}-k_0\mu]^2}{2k_0\mu
+\frac{2}{3}[(\ln n)^{\frac{1}{3}}-k_0\mu]}\right\}\nonumber\\
\!\!\!\!\!& = &\!\!\!\!\!0 \quad \mbox{a.a.s.}
\end{eqnarray}

Finally, we show fact (ii).  Suppose that for each $i=1,2,...$,
there exists $n_i$ and $v\in G(\mathcal{X}^{(d)}_{n_i};1)
\in{\mathcal T}_{n_i}$ with $n_i\geq
|\Gamma_v|>\frac{3-2\epsilon}{\epsilon^2}\ln n_i$ and a realization
of the ``active-saturated" process starting at $v$ such that
$|\Gamma_v|-1\geq k'\geq i$. As $i\rightarrow \infty$,
$k'\rightarrow \infty$, $|\Gamma_v| \rightarrow \infty$, and
$n_i\rightarrow \infty$.  However, this holds with probability 0.
This completes our proof. \qed


Note that by either ignoring border effects or by taking
$n\rightarrow \infty$, $C_t^{(d)}$ (defined for
$G(\mathcal{X}^{(d)}_n;1)$) is also the $t$-th order cluster
coefficient for the infinite Poisson random geometric graph
$G(\mathcal{H}^{(d)}_{\lambda};1)$.  Thus by Theorem
\ref{Theorem-RGG-Lower-Bound}, we obtain the following important
corollaries giving new improved lower bounds on the critical density
for $d$-dimensional Poisson random geometric graphs.

\vspace{0.1in}%
\begin{corollary}\label{Corollary-RGG-Critical-Mean-Degree}
Let $\mu_c^{(d)}$ and $\lambda^{(d)}_c$ be the critical mean degree
and critical density for $G(\mathcal{H}^{(d)}_{\lambda};1)$,
respectively, where $d \geq 2$. Then for all $t \geq 3$,
\begin{equation}
\mu_c^{(d)} \geq \frac{1}{1-C_t^{(d)}}, \quad\mbox{and}\quad \lambda^{(d)}_c\geq \frac{1}{V^{(d)}
\left(1-C_t^{(d)}\right)}. \label{lowerbound}
\end{equation}
\end{corollary}
\vspace{0.1in}%

\emph{Proof}:  Follows immediately from Theorem
\ref{Theorem-RGG-Lower-Bound} and~\eqref{mu-lambda}. \qed

In particular, for two-dimensional Poisson random geometric graphs, substituting $V^{(2)}=\pi$ and
$C^{(2)}_3=1-\frac{3\sqrt{3}}{4\pi}=0.5865...$ into~(\ref{lowerbound}), we have the following
corollary.

\vspace{0.1in}%
\begin{corollary}\label{Corollary-RGG-Critical-Density}
The critical mean degree $\mu_c^{(2)}$ and the critical density $\lambda^{(2)}_c$ for
$G(\mathcal{H}^{(2)}_{\lambda};1)$ satisfy $\mu_c^{(2)} \geq 2.419...$ and $\lambda_c^{(2)}\geq
0.7698...$.
\end{corollary}
\vspace{0.1in}%

Note that if we use high-order cluster coefficients $C_t^{(2)},
t\geq 4$, computed by numerical methods, e.g., $C^{(2)}_4 \approx
0.6012$ and $C^{(2)}_5 \approx 0.6179$, we can obtain further
improved (approximate) lower bounds: $\mu_c \gtrsim 2.617$, and
$\lambda_c\gtrsim 0.883$.\footnote{Of course, if techniques are
developed to compute the higher-order cluster coefficients in closed
form, we would automatically obtain even tighter analytical lower
bounds.}

By applying $C_3^{(3)}=0.4688...$ given by (\ref{C-3-3}), we have

\vspace{0.1in}%
\begin{corollary}\label{Corollary-RGG-3-Dimension}
The critical mean degree $\mu_c^{(3)}$ and the critical density $\lambda^{(3)}_c$ for
$G(\mathcal{H}_{\lambda}^{(3)};1)$ satisfy $\mu_c^{(3)} \geq 1.412...$ and $\lambda_c^{(3)}\geq
0.4494...$.
\end{corollary}
\vspace{0.1in}%

This lower bound is close to the known results obtained by simulation---0.65 \cite{AlDrBa90}, and
it is the best known analytical lower bound on $\lambda_c^{(3)}$.

\section{Conclusion}

We have established a new class of analytical lower bounds on the
critical density $\lambda_c^{(d)}$ for percolation in
$d$-dimensional Poisson random geometric graphs.  These analytical
lower bounds are the tightest known to date, and reveal a deep
underlying relationship between the cluster coefficient and the
critical density in continuum percolation.

\section*{Acknowledgement}

We would like to thank the anonymous referees for their helpful comments and suggestions, which
improved the paper.

This research is supported in part by Army Research Office (ARO) grant W911NF-06-1-0403.

\bibliography{PercolationTopic2}

\begin{thebibliography}{10}

\bibitem{BoNrFrMe03}
L.~Booth, J.~Bruck, M.~Franceschetti, and R.~Meester, ``Covering algorithms,
  continuum percolation and the geometry of wireless networks,'' {\em Annals of
  Applied Probability}, vol.~13, pp.~722--741, May 2003.

\bibitem{FrBoCoBrMe05}
M.~Franceschetti, L.~Booth, M.~Cook, J.~Bruck, and R.~Meester, ``Continuum
  percolation with unreliable and spread out connections,'' {\em Journal of
  Statistical Physics}, vol.~118, pp.~721--734, Feb. 2005.

\bibitem{DoFrTh05}
O.~Dousse, M.~Franceschetti, and P.~Thiran, ``Information theoretic bounds on
  the throughput scaling of wireless relay networks,'' in {\em Proc. IEEE
  INFOCOM'05}, Mar. 2005.

\bibitem{DoBaTh05}
O.~Dousse, F.~Baccelli, and P.~Thiran, ``Impact of interferences on
  connectivity in ad hoc networks,'' {\em IEEE Trans. Network.}, vol.~13,
  pp.~425--436, April 2005.

\bibitem{Gi61}
E.~N. Gilbert, ``Random plane networks,'' {\em J. Soc. Indust. Appl. Math.},
  vol.~9, pp.~533--543, 1961.

\bibitem{Ha85}
P.~Hall, ``On continuum percolation,'' {\em Annals of Prob.}, vol.~13,
  pp.~1250--1266, 1985.

\bibitem{MeRo96}
R.~Meester and R.~Roy, {\em Continuum Percolation}.
\newblock New York: Cambridge University Press, 1996.

\bibitem{Pe03}
M.~Penrose, {\em Random Geometric Graphs}.
\newblock New York: Oxford University Press, 2003.

\bibitem{Gr99}
G.~Grimmett, {\em Percolation}.
\newblock New York: Springer, second~ed., 1999.

\bibitem{BoRi06}
B.~Bollob{\'{a}}s and O.~Riordan, {\em Percolation}.
\newblock New York: Cambridge University Press, 2006.

\bibitem{QuToZi00}
J.~Quintanilla, S.~Torquato, and R.~M. Ziff, ``Efficient measurement of the
  percoaltion threshold for fully penetrable discs,'' {\em Physics A}, vol.~86,
  pp.~399--407, 2000.

\bibitem{BaBoWa05}
P.~Balister, B.~Bollob{\'{a}}s, and M.~Walters, ``Continuum percolation with
  steps in the square or the disc,'' {\em Random Structures Algorithms},
  vol.~26, pp.~392--403, 2005.

\bibitem{JaLuRu00}
S.~Janson, T.~Luczak, and A.~Ruci\'{n}ski, {\em Random Graphs}.
\newblock New York: John Wiley \& Sons, 2000.

\bibitem{Fe57}
W.~Feller, {\em An Introduction to Probability Theory and its Applications},
  vol.~1.
\newblock New York: John Wiley \& Sons, 1957.

\bibitem{DaCh02}
J.~Dall and M.~Christensen, ``Random geometric graphs,'' {\em Phy. Rev. E},
  vol.~66, no.~016121, 2002.

\bibitem{AlDrBa90}
U.~Alon, A.~Drory, and I.~Balberg, ``Systematic derivationof percolation
  thresholds in continuum systems,'' {\em Physical Review A}, vol.~42,
  pp.~4634--4638, 1990.

\end{thebibliography}
\bibliographystyle{ieeetr}

\section*{Appendix I}

\begin{proposition}\label{Proposition-Coupling}
Suppose random variables $X_i, i=1,...,m$ satisfy the following conditions: (i)
$\mbox{Pr}\{X_1\geq x\}\leq \mbox{Pr}\{X_1^+\geq x\},\;\forall x$; (ii)
$\mbox{Pr}\{X_l\geq x|X_1=x_1,...,X_{l-1}=x_{l-1}\}\leq \mbox{Pr}\{X_l^+\geq
x\},\;\forall x,x_1,...,x_{l-1}$; (iii) $X_i^+, i=1,...,m$ are independent of each other
and of $X_i, i=1,...,m$. Then
\[
\mbox{Pr}\left\{\sum_{i=1}^mX_i\geq z\right\}\leq \mbox{Pr}\left\{\sum_{i=1}^m X_i^+\geq
z\right\}, \;\forall z
\]
\end{proposition}
\vspace{0.1in}%

\emph{Proof:} It suffices to show the result for $m=2$. Since $\mbox{Pr}\{X_2\geq
y|X_1=x_1\}\leq\mbox{Pr}\{X_2^+\geq y|X_1=x_1\},\; \forall x_1,y$, $\mbox{Pr}\{X_1+X_2\geq
z|X_1=x_1\}\leq\mbox{Pr}\{X_1+X_2^+\geq z|X_1=x_1\},\; \forall x_1,z$, and thus
\[
\begin{array}{lll}
\vspace{+.1in}
& &\mbox{Pr}\{X_1+X_2\geq z\} \\
\vspace{+.1in}
& =  & \sum_{x_1}\mbox{Pr}\{X_1=x_1\}\mbox{Pr}\{X_1+X_2\geq z|X_1=x_1\}\\
\vspace{+.1in}
&  \leq  & \sum_{x_1}\mbox{Pr}\{X_1=x_1\}\mbox{Pr}\{X_1+X_2^+\geq z|X_1=x_1\} \\
\vspace{+.1in}
&  =  & \mbox{Pr}\{X_1+X_2^+\geq z\}\\
\vspace{+.1in}
&  =  &\sum_{x_2^+}\mbox{Pr}\{X_2^+=x_2^+\}\mbox{Pr}\{X_1+x_2^+\geq z|X_2^+=x_2^+\}\\
\vspace{+.1in}
&  =  & \sum_{x_2^+}\mbox{Pr}\{X_2^+=x_2^+\}\mbox{Pr}\{X_1\geq z-x_2^+\}\\
\vspace{+.1in}
&  \leq & \sum_{x_2^+}\mbox{Pr}\{X_2^+=x_2^+\}\mbox{Pr}\{X_1^+\geq z-x_2^+\}\\
\vspace{+.1in}
&  =  & \sum_{x_2^+}\mbox{Pr}\{X_2^+=x_2^+\}\mbox{Pr}\{X_1^++X_2^+\geq z|X_2^+=x_2^+\}\\
&  =  & \mbox{Pr}\{X_1^++X_2^+\geq z\}.
\end{array}
\]
\qed

\end{document}